\newcount\notenumber

\def\note{\advance\notenumber by 1
\footnote{$^{(\the\notenumber)}$}}

\def\Q{{\bf Q}}

\def\N{{\bf N}}
\def\Z{{\bf Z}}

\def\x{{\bf x}}

\def\V{{\cal V}}

\def\C{{\bf C}}
\def\G{{\bf G}}

\def\Z{{\bf Z}}


\font\title=cmr10 scaled 1200

\centerline{\title On the greatest prime factor of
$(ab+1)(ac+1)$}\bigskip

\centerline{P. Corvaja\hskip 2truecm U. Zannier}\vskip
1truecm

\noindent{\bf Abstract.} We prove  that for integers $a>b>c>0$, the
greatest prime factor of $(ab+1)(ac+1)$ tends to infinity with $a$.
 In particular, this settles a  conjecture raised by 
Gy\"ory, Sarkozy and Stewart, predicting the same conclusion for 
the product
$(ab+1)(ac+1)(bc+1)$.
\bigskip

In the paper [GSS], G\"yory, Sarkozy and Stewart conjectured that, for
positive integers $a>b>c$, the greatest prime factor of the
product $(ab+1)(ac+1)(bc+1)$ tends to infinity as
$a\rightarrow\infty$. In the direction of
this conjecture, some partial  results and analogues for more than three
integers were obtained in the papers  [GS] and [ST]; for example, in this
last paper the conclusion is proved (through the use of Baker's method, and
thus in en ``effective" way)  under the assumption that $\log c/\log
a\rightarrow 0$. The paper [GS] instead proved the special case when some
number in the set $a,b,c,a/b,b/c,c/a$ is an $S$-unit, and quantitative
improvements on [GS] come from a paper by Bugeaud [B].  However the general
case remained open. \medskip

In this note we show how the methods introduced in [CZ] and [BCZ]  may be
useful for such problems. In fact, we shall  give an affirmative answer
actually to a strengthening of the above conjecture. Namely, we prove that the
greatest prime factor of $(ab+1)(ac+1)$ tends to infinity with $a$.  

Actually, it will be clear from the proof that
arguments  similar to ours in fact lead to much more general conclusions (see
also Remark 1 at the end). However here we shall focus only on the basic
case just  mentioned. 

We shall prove the above statement (and so the conjecture) in the following
equivalent form:\medskip

\noindent{\bf Theorem.}\enspace{\it Let $S$ be a finite set of prime numbers.
Then there exist only finitely many  triples of positive integers $a>b>c$ such
that the product $(ab+1)(ac+1)$ has all of its prime factors in
$S$.}\medskip

As mentioned above, our proof will use the Subspace Theorem, leading 
(contrary to [ST]) to an 
ineffective result, in the sense that it will not  provide an explicit lower
bound for the greatest prime factor in question (i.e. our proof  will not
provide an explicit upper bound for $a$ in terms of $S$). However, for a given
$S$  it  would be possible to obtain an explicit   upper bound for the number
of triples in the above statement. \medskip

\noindent{\it Acknowledgement.} We thank Yann Bugeaud for drawing our
attention to the problem and for providing us  with some relevant
references.\medskip

\noindent{\it Proof of theorem.} We shall use the same letter $S$ as in the
statement to mean also the set of valuations of $\Q$ associated to  the primes
in $S$. We also suppose that $S$ includes the infinite valuation.

 Our arguments depend on the Subspace Theorem, as well as on   a result by
Liardet, appealed to in the final part.  (This is  essentially a particular
case of the former ``Lang conjecture for algebraic tori", proved later by M.
Laurent in full generality.) For the reader's convenience we recall at once a 
formulation  of the Subspace Theorem relevant to us. 
(For a proof see [S1,2].)\medskip

\noindent{\bf Subspace Theorem.} {\it Let $S$ be a finite set of absolute
values of $\Q$, including $\infty$ (normalized so that $|p|_p=p^{-1}$) and let
$N\in\N$. For $w\in S$,  let $L_{1w},\ldots ,L_{Nw}$ be linearly independent
linear forms in $N$ variables, with rational coefficients, and let $\delta
>0$. Then the solutions $\x =(x_1,\ldots ,x_N)\in\Z^{N}$ of the inequality
$$
\prod_{w\in S}\prod_{j=1}^N|L_{jw}(\x)|_w< (\max |x_i|)^{-\delta} 
$$
are contained in finitely many proper subspaces of $\Q^{N}$.}\medskip

Now with the proof. We argue by contradiction and we let $\Sigma$ denote an
infinite set of triples of positive integers $(a,b,c)$ with $a>b>c$ such
that   the two integers $u:=ab+1$, $v:=ac+1$  are
$S$-units (namely they are entirely composed of primes from $S$). \medskip

Further, we put 
$$
y_1:={u-1\over v-1}={b\over c},\qquad y_2={u^2-1\over v-1}={(u+1)b\over
c}, 
$$
so $y_1,y_2$ are rational numbers with denominator at most $c$. Similarly to
[CZ, Thm. 1] and [BCZ], we shall now approximate $y_1,y_2$ with suitable
linear combinations of $S$-units, which will allow an application of the
Subspace Theorem. 

First,  observe the approximation\note{where for our purpose the number $5$ 
could be replaced by any larger integer} 
$$
{1\over v-1}={1\over v(1-v^{-1})}=\sum_{n=1}^\infty
v^{-n}=\sum_{n=1}^5v^{-n}+O(v^{-6}).
$$
On multiplying by $u^j-1$, for $j=1,2$, we thus obtain
$$
\left|y_j+\sum_{n=1}^5v^{-n}-\sum_{n=1}^5u^jv^{-n}\right|\ll
u^jv^{-6},\qquad j=1,2,
$$
where the implied constant is absolute. In turn, this is equivalent to
$$
\left|v^5y_j+\sum_{n=1}^5v^{5-n}-\sum_{n=1}^5u^jv^{5-n}\right|\ll
u^jv^{-1},\qquad j=1,2.\eqno(1)
$$
Let $\sigma_1,\ldots ,\sigma_{15}$ denote the integers $u^jv^{5-n}$ for
$j=0,1,2$ and $n=1,\ldots ,5$, in some order. Naturally, these integers depend
on $(a,b,c)\in\Sigma$.  Then (1) may be rewritten in the form
$$
\left|v^5y_j+\sum_{i=1}^{15}\alpha_{ji}\sigma_i\right|\ll u^jv^{-1}\qquad
j=1,2,\eqno(2) 
$$
for suitable rational numbers $\alpha_{ji}$.\medskip

We now   introduce linear forms $L_{jw}$ in $17$ 
variables $Y_1,Y_2,X_1\ldots ,X_{15}$, for $j=1,\ldots ,17$ and $w\in S$. We
set  
$$
L_{j\infty}=Y_j+\sum_{i=1}^{15}\alpha_{ji}X_i,\qquad
L_{jw}=Y_j\quad {\rm for}\ w\neq\infty ,\qquad j=1,2, 
$$
while for $j=3,\ldots ,17$ and any $w\in S$ we put simply
$L_{jw}=X_{j-2}$. 

Plainly, for each $w\in S$ the linear forms $L_{jw}$, $j=1,\ldots ,17$ are
linearly independent.\medskip

We also define, for each $(a,b,c)\in\Sigma$, a vector 
$\x=(x_1,\ldots ,x_{17})\in\Z^{17}$ by
$$
 \x=(x_1,\ldots ,x_{17})=(cv^5y_1,cv^5y_2,c\sigma_1,\ldots ,c\sigma_{15}).
$$
Observe that the $x_i$ so defined are in fact integers. Inequalities (2)
translate into
$$
|L_{j\infty}(\x)|_\infty\ll cu^jv^{-1},\qquad j=1,2.\eqno(3).
$$
Also, for $j=1,2$ we have
$$
\prod_{w\in S\setminus\infty}|L_{jw}(\x)|_w=
\prod_{w\in S\setminus\infty}|cy_jv^5|_w\le 
\prod_{w\in S\setminus\infty}|v^5|_w=v^{-5},
$$
the  inequality holding because $cy_j$ is an integer and the last equality
following from the product formula, since $v$ is a positive  $S$-unit.
Combining with (3) yields
$$
\prod_{w\in S}|L_{jw}(\x)|_w\ll cu^jv^{-6},\qquad j=1,2.\eqno(4)
$$
On the other hand, if $j=3,\ldots ,17$ we have (by the product formula again)
$$
\prod_{w\in S}|L_{jw}(\x)|_w=\prod_{w\in S}|c\sigma_{j-2}|_w\le c,\eqno(5)
$$
since $c$ is an integer and since $\sigma_1,\ldots ,\sigma_{15}$ are $S$-units.
Combining (4) and (5) yields
$$
\prod_{j=1}^{17}\prod_{w\in S}|L_{jw}(\x)|_w\ll c^{17}u^3v^{-12}.\eqno(6)
$$
Now, we have   $u=ab+1\le a^2$ and $v=ac+1> ac$, whence 
$c^{17}u^3v^{-12}\le c^{5}a^{-6}< 
a^{-1}$. Also, $\max |x_i|\le u^2v^5c\le a^{15}$, so (6) gives
$$
\prod_{i=1}^{17}\prod_{w\in S}|L_{iw}(\x)|_w\ll (\max |x_i|)^{-{1\over 15}}.
$$
Therefore we may apply the Subspace Theorem (say with $\delta =1/16$) and
deduce that the vectors $\x$ in question all lie on the union of finitely many
rational proper linear  subspaces of $\Q^{17}$. In particular, we may assume
that there exist rationals $\eta_1,\eta_2,\gamma_1,\ldots ,\gamma_{15}$,  
not all zero and  such that  for infinitely many triples $(a,b,c)\in\Sigma$,
we have   
$$
\eta_1y_1v^5+\eta_2y_2v^5+\sum_{i=1}^{15}\gamma_i\sigma_i=0.
$$
Recalling the definition of $y_i$ and
$\sigma_i$ we  derive an equation 
$$
\eta_1{v^5(u-1)\over v-1}+\eta_2{v^5(u^2-1)\over v-1}
+\sum_{j,n}\rho_{jn}u^jv^{5-n}=0,\eqno(7)
$$
valid for an infinity of triples in $\Sigma$, where the summation is
over $j=0,1,2$ and $n=1,\ldots ,5$ and where the coefficients
$\eta_1,\eta_2,\rho_{jn}$ are fixed rationals, not all zero. \medskip

Consider now the algebraic curve $\V$   (not necessarily irreducible) defined
in $\G_m^2$ by the equation 
$$
\eta_1(U-1)+\eta_2(U^2-1)
+(V-1)\sum_{j,n}\rho_{jn}U^jV^{-n}=0,
$$
where $U,V$ are  variables. We pause to show that this equation is
nontrivial, so it defines in fact a curve. Suppose on the contrary that  the
left side  vanishes identically.  Then $\eta_1(U-1)+\eta_2(U^2-1)$ would be
divisible by $V-1$, (which is coprime with the denominators of all the other
terms). But this would imply $\eta_1=\eta_2=0$, yielding 
$\sum_{j,n}\rho_{jn}U^jV^{-n}=0$. Now, looking at nonzero terms with maximal
$n$, we conclude that $\rho_{jn}=0$ for all $j,n$ in question, a
contradiction.\medskip

In view of (7), $\V$ contains infinitely many
points $(u,v)$ in the finitely generated group of points whose coordinates
are $S$-units in $\Q$. By a theorem of Liardet (see e.g. 
[L, Thm. 7.3, p. 207]) all but finitely many of such points lie in a certain 
finite  union of  translates of algebraic subtori of $\G_m^2$.
So, we may assume that some such translate contains an infinity of the points
in question.  But it is an easy well-known fact that such a translate is
defined by some equation $U^pV^q=h$, where $p,q$ are coprime integers not both
zero and where $h\in\C^*$ (see also the quoted theorem in [L]). Therefore we
would have $u^pv^q=h$ for some fixed $h$ and  an infinity of triples in
$\Sigma$. This plainly implies that $h$ is rational and that $pq <0$. So,
writing $h=h_1/h_2$  and replacing $p,q$ by $\pm p, \pm
q$, we may assume that $h_1u^p=h_2v^q$, where now $p,q, h_1,h_2$ are
certain fixed positive integers with $(p,q)=1$. But $u\equiv v\equiv 1\pmod a$,
whence $h_1\equiv h_2\pmod a$ for an infinity of integers $a$. This implies
$h_1=h_2$, so the equation is in fact $u^p=v^q$. This  implies $u=t^q$,
$v=t^p$ for  a suitable integer $t$,
depending on $a,b,c$. Now, the GCD$(t^p-1,t^q-1)$ is bounded by a constant
multiple of $t-1$, since the polynomials  ${x^p-1\over x-1}$ and ${x^q-1\over
x-1}$ are coprime. On the other hand the GCD$(u-1,v-1)$ is a multiple of $a$,
whence  
$$ 
a\ll t-1\le u^{1/q}\le a^{2/q}.
$$
Therefore $q\le 2$, which forces $p=1, q=2$, i.e. $u=v^2>a^2$.  This would
imply however $b\ge a$, a contradiction which proves the theorem.///\medskip

\noindent{\bf Remark 1.} A combination of the above arguments with the 
more refined technique
of [BCZ] leads  to an improvement on that paper. In such a  sharpening one
can   show the following: {\it Let $\epsilon >0$ and let  $S$ be a finite set
of primes. Then, for  pairs $(u,v)$ of multiplicatively independent $S$-units,
we have   {\rm GCD}$(u-1,v-1)\ll_{\epsilon ,S} (\max (u,v))^{\epsilon}$}. (The
result of [BCZ] is  the special case  $u=a^n$, $v=b^n$, where $a,b$ are fixed
multiplicatively independent integers and $n$ varies through $\N$.) \medskip

\noindent{\bf Remark 2.} A somewhat different 
solution to the original conjecture is  as follows. Write $ab+1=r$, $ac+1=s$,
$bc+1=t$, so $r,s,t$ are $S$-units. We find $(abc)^2=rst-rs-rt-st+r+s+t-1$.
This equation may be treated as in [CZ], by expanding with the binomial 
theorem the square root of the right side. This works {\it provided} the term
$rst$ is ``dominant", which here means that $rs<(rst)^{1-\delta}$, for some
fixed positive $\delta$. In turn, it suffices that $b$ is larger than a fixed
positive power of $a$. On the other hand, if this is not the case, then
$\log c/\log a$ tends to zero and the result of [ST] applies.\medskip

\bigskip

{\bf References}

\medskip

\item{[B]} Y. Bugeaud, On the greatest prime factor of $(ab+1)(ac+1)(bc+1)$,
{\it Acta Arith.}, {\bf 86} (1998), 45-49.\smallskip

\item{[BCZ]} Y. Bugeaud, P. Corvaja, U. Zannier - An upper bound for  the
G.C.D. of $a^n-1$ and $b^n-1$, preprint, 2001.\smallskip

\item{[CZ]} P. Corvaja, U. Zannier - Diophantine equations with power sums 
and universal Hilbert
  sets, {\it Indagationes Math.}, {\bf 9} (1998), 317-332. \smallskip

\item{[GS]} K. Gy\"ory, A. Sarkozy - On prime factors of integers of the form 
$(ab+1)(ac+1)(bc+1)$, {\it Acta Arith.}, {\bf 79} (1997).\smallskip

\item{[GSS]} K. Gy\"ory, A. Sarkozy, C.L. Stewart - On the number of
prime factors of integers of the form $ab+1$, {\it Acta Arith.}, {\bf
74} (1996), 365-385.\smallskip

\item{[L]} S.Lang - {\it Fundamentals of Diophantine Geometry}, 
Springer-Verlag, 1983.\smallskip

\item{[S1]} W.M. Schmidt - {\it Diophantine Approximation}, Springer-Verlag
LNM {\bf 785}, 1980. \smallskip

\item{[S2]} W.M. Schmidt - {\it Diophantine Approximations 
and Diophantine Equations}, Springer-Verlag
LNM {\bf 1467}, 1991.\smallskip

\item{[ST]} C.L. Stewart, R. Tijdeman - On the greatest prime factor of
$(ab+1)(ac+1)(bc+1)$, {\it Acta Arith.}, {\bf 79} (1997).\smallskip

\vfill

 {Umberto Zannier} \hfill Pietro Corvaja 

 {I.U.A.V. - DCA} \hfill Dipartimento di Mat. e Inf.  

{S. Croce, 191} \hfill via delle Scienze, 206 

 {30135 VENEZIA (ITALY)} \hfill 33100 UDINE (ITALY) 

  {zannier@iuav.it}\hfill corvaja@dimi.uniud.it

\end